\documentclass[11pt]{article}
\usepackage{graphicx}

\def\ds{\displaystyle}

\def\TV{\hbox{Tot.Var.}}
\def\L{{\bf L}}

\def\ve{\varepsilon}

\def\vsk{\vskip 4em}
\def\v{\vskip 1em}
\def\O{{\cal O}}

\def\bega{\begin{array}}
\def\enda{\end{array}}
\def\begi{\begin{itemize}}
\def\endi{\end{itemize}}

\def\bel{\begin{equation}\label}
\def\eeq{\end{equation}}
\def\sqr#1#2{\vbox{\hrule height .#2pt
\hbox{\vrule width .#2pt height #1pt \kern #1pt
\vrule width .#2pt}\hrule height .#2pt }}
\def\square{\sqr74}
\def\endproof{\hphantom{MM}\hfill\llap{$\square$}\goodbreak}

\begin{document}
\title{\bf Lack of BV Bounds for Approximate Solutions
to the p-system with Large Data}
\vsk

\author{Alberto Bressan$^{(*)}$, Geng Chen$^{(**)}$, and Qingtian Zhang$^{(*)}$\\    \\
(*) Department of Mathematics, Penn State University
University Park, Pa.~16802, U.S.A.\\
(**) School of Mathematics
Georgia Institute of Technology
Atlanta, Ga.~30332, U.S.A..
\\
e-mails:~ bressan@math.psu.edu~,~gchen73@math.gatech.edu~, ~zhang\_q@math.psu.edu}

\maketitle

\begin{abstract} We consider front tracking approximate solutions to the p-system of
isentropic gas dynamics. At interaction times,  the outgoing wave fronts
have  the same strength as in the exact solution of the Riemann problem,
but some error is allowed in their speed. For large BV initial data,
we construct examples showing that the total variation of these approximate solutions
can become arbitrarily large, or even blow up in finite time.
This happens even if
the density of the gas remains uniformly positive.
\end{abstract}

\section{Introduction}
\label{sec:0}
\setcounter{equation}{0}

For hyperbolic systems of conservation laws in one space dimension,
a satisfactory
existence-uniqueness theory is currently available for entropy weak solutions
with small total variation \cite{Bbook, BLY}.
The well-posedness of the Cauchy problem
holds also in the case of large data, as long as the total variation remains
bounded \cite{BC1, Lw}.  A major remaining open problem is whether, for large
BV initial data,
the total variation remains uniformly
bounded or can blow up in finite time.
Examples of finite time blow up have been constructed in \cite{BJ, J}.  However, these
systems do not come from physical models and do not admit a strictly convex entropy.

For initial data having small total variation,
regardless of the order in which
wave fronts cross each other, the Glimm interaction estimates \cite{G} show that
the total strength of waves remains small for all times.  There are few examples
of hyperbolic systems where uniform BV estimates hold also
for solutions with large data \cite{N, T}.
In the present paper
we study BV bounds for the p-system of isentropic gas dynamics in Lagrangean variables:
\bel{1}\left\{\begin{array}{rl} v_t - u_x&=~0\,,\cr
u_t + p(v)_x&=~0\,.\end{array}\right.\eeq
Here $u$ is the velocity, $\rho$ is the density,  $v=\rho^{-1}$ is specific volume,
while $p=p(v)$ is the pressure.
Our main concern is whether, for front tracking approximate solutions to
the p-system, uniform BV estimates  can be established.  More precisely,
we study the
following question:

\begi
\item[{\bf (Q)}] {\it 
Consider a front tracking approximation for (\ref{1}) with large BV initial data.
Assume that, at each interaction time,  the outgoing wave fronts
have  the same strength as in the exact solution of the Riemann problem
but some error is allowed in their speed.
Can the total variation
of such approximate solution become arbitrarily large?}
\endi

In  this paper, some examples will be constructed, showing
that  the total strength of waves in a front tracking approximation
can indeed approach infinity in finite time.   This confirms the non-existence
of a Lyapunov functional which decreases at every wave-front interaction,
as proved in \cite{CJ}.

For interactions occurring near vacuum, it was already noticed in
\cite{LS} that uniform Glimm-type estimates were no longer valid.
It thus comes as no surprise that an approximate front tracking solution
can be constructed, where the total strength of waves (measured by the
change in Riemann invariants) blows up in finite time.
Remarkably, our last two examples show that
an arbitrarily large amplification of the total
variation is still possible even if the gas density remains uniformly positive.

It should be clear that the present counterexamples do not prove that,
for large BV solutions  of the $p$-system, the total variation can
blow up in finite time.   Indeed, we still conjecture that global BV bounds
do hold.    Our analysis only shows that such BV bounds cannot be proved
by wave interaction estimates alone, and additional properties of solutions
must be taken into account.  Apparently, the decay of rarefaction waves
due to genuine nonlinearity \cite{Bbook, BC2,BY, GL} should be used in a crucial way.
In the last section we revisit two of the earlier examples and show that,
if such decay is taken into account,  these specific interaction patterns
do not produce blow up.

\section{Wave interactions for the p-system}
\label{sec:1}
\setcounter{equation}{0}

We review here some standard properties of characteristic curves and of shock curves.
For details we refer to \cite{CJ,Sm}.
To simplify the computations,
 we assume that in (\ref{1})  the pressure has the special form
$$p(v)~=~{1\over 3v^3}~= ~{\rho^3\over 3}\,,\qquad\qquad p'(v) ~=~-{1\over v^4}\,.$$
Smooth solutions of (\ref{1}) satisfy the quasilinear system
\bel{ql}
\left\{ \bega{rl} \rho_t + \rho^2 u_x &=~0,\cr\cr
u_t + \rho^2 \rho_x &=~0,\enda\right.\eeq
with characteristic speeds
\bel{2}\lambda~=~\pm\rho^2\,, \eeq
The variables
\bel{RI}
w_1~=~\rho-u\,,\qquad\qquad w_2~=~\rho + u\,,\eeq
provide a coordinate system of Riemann invariants, in the $u$-$\rho$ plane.
\begin{figure}[htbp]
\centering
  \includegraphics[scale=0.4]{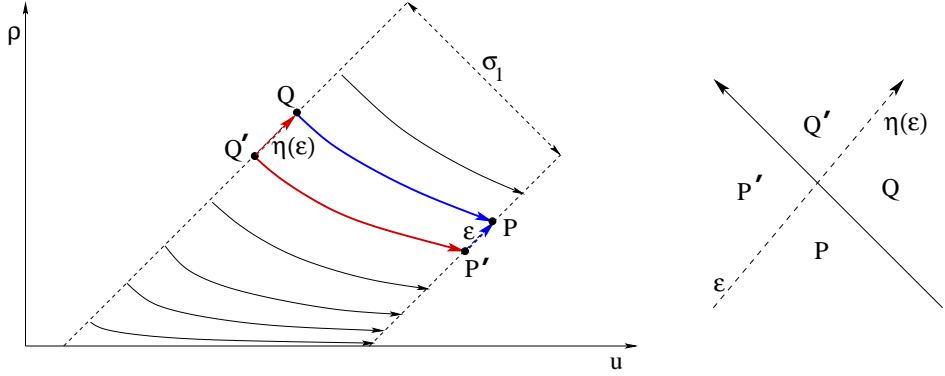}
    \caption{{\small Interaction of a 1-shock and a 2-rarefaction}}
\label{f:hyp81}
\end{figure}

A shock with left state $(u_-,\rho_- )$ and right state $(u_+,\rho_+)$,
traveling with speed
$\lambda$,
satisfies the Rankine-Hugoniot equations
$$\left\{ \begin{array}{rl} \ds \lambda\left({1\over \rho_+}- {1\over \rho_-}
\right)&=~u_--u_+\,,\cr
&\cr
\lambda(u_+-u_-) & = ~p(v_+) - p(v_-) ~=\ds
~{\rho_+^3\over 3} - {\rho_-^3\over 3}\,.\end{array}
\right.$$
Hence
\bel{RH2}u_+-u_- ~=~- \sqrt{ \left({1\over\rho_+}- {1\over \rho_-}\right)\left(
{\rho_-^3\over 3}-{\rho_+^3\over 3}\right)}\,,\eeq
\bel{RH5}\lambda ~=~\pm {1\over\sqrt 3}
{\ds\sqrt{\ds {\rho^3_-}-{\rho^3_+} \over \ds
{1\over\rho_+}- {1\over \rho_-}}}\,.\eeq
Setting
$$ s~=~u_--u_+\,,\qquad\qquad \theta~=~{\rho_+\over\rho_-}\,,$$
from (\ref{RH2}) it follows
\bel{RH4}
s~=~\rho_-\, \sqrt{(1-\theta)(1-\theta^3)\over 3\theta}\,.
\eeq

For the p-system, the interaction of wave fronts has been thoroughly studied
\cite{Sm, CJ}. For reader's convenience, we derive here
two elementary estimates which will be used in the sequel.

Consider a 1-shock with left and right states (see Fig.~\ref{f:hyp81})
$$P~=~(u_-, \rho_-),\qquad\qquad Q~=~(u_+, \rho_+)\,,$$
with strength
$$\sigma_1~=~|(\rho_+-u_+) - (\rho_--u_-)|$$ measured by the change in the
1-Riemann invariant.  Assume that this shock
crosses a small 2-wave (compression or rarefaction)
of strength $\sigma_2=\ve$.
We seek an estimate on the size of the outgoing waves, up to leading order.
Let
$$P'~=~(u_--\ve, \,\rho_--\ve),\qquad\qquad Q'~=~\Big(u_+-\eta(\ve),
\,\rho_+-\eta(\ve)\Big)
$$
be the left and right states across the 1-shock after the interaction.
Setting
\bel{ff}s(\ve)~=~s-\ve+\eta(\ve),\qquad\qquad\theta(\ve)~=~{\theta\rho_- -\eta(\ve)
\over \rho_--\ve}\,,\qquad\qquad \psi(\theta)~=~
{(1-\theta)(1-\theta^3)\over 3\theta}\eeq
from (\ref{RH4}) it follows
\bel{RH7}
s(\ve)~=~(\rho_--\ve)\, \sqrt{\psi(\theta(\ve))}\,.
\eeq
Differentiating w.r.t.~$\ve$, at $\ve=0$ we obtain
\bel{der}\bega{rl}
\eta'-1&\ds=~-\sqrt{\psi(\theta)} + {\rho_-\over 2\sqrt{\psi(\theta)}}\cdot \psi'(\theta)\theta'
\cr\cr
&=\ds-{s\over\rho_-} + {\rho_-^2\over 2s}\cdot {-1-2\theta^3+3\theta^4\over 3\theta^2}
\cdot {\theta-\eta'\over \rho_-}\,.\enda\eeq
Two cases are relevant to our analysis.
\v
CASE 1: The 1-shock has small amplitude.   In this case,
since shock and rarefaction curves coincide up to second order,
for $s$ small we have the expansion
\bel{thex}
\theta~=~1+{s\over \rho_-}+\O(s^3)\,.\eeq
Inserting (\ref{thex}) in (\ref{der}) we obtain
\bel{epex}
\eta'~=~1+\frac{1}{3} \left(\frac s{\rho_-}\right)^3 + \O(s^4)\eeq
\v
CASE 2: The 1-shock has a fixed strength $\sigma_1$, while the density $\rho_-$
approaches zero.
By definition, the strength is computed by
\bel{ss1}\sigma_1~=~(\theta-1)\rho_-+s~=~\left[ (\theta-1)+\sqrt{\theta^3-\theta^2-1+
\theta^{-1}\over 3}\right]\rho_-\,.\eeq
As $\rho_-\to 0$, from (\ref{ss1}) it follows that $\theta\to \infty$.
Indeed, dropping lower order terms we find
\bel{vac}{1\over\sqrt 3}\theta^{3/2}\rho_-~\approx~\sigma_1\,,
\qquad\qquad
\theta~\approx~\left(\sqrt 3 \sigma_1\over \rho_-\right)^{2/3}\,.
\eeq
Inserting (\ref{vac}) in (\ref{der}) and retaining only leading order terms,
we obtain
$$\eta'-1~\approx~{-\sigma_1\over \rho_-} + {\rho_-\over 2\sigma_1}\cdot \theta^2\,
(\theta-\eta')\,,$$
$$\eta'\approx~\frac{\ds
1-\frac{\sigma_1}{\rho_-}+\frac32\frac{\sigma_1}{\rho_-}}{\ds
1+\frac{\rho_-}
{2\sigma_1}\Big(\sqrt 3\frac{\sigma_1}{\rho_-}\Big)^{4/3}}~\approx~
\frac1{3^{2/3}}\Big(\frac{\sigma_1}{\rho_-}\Big)^{2/3},$$
\bel{ep}
\eta'
~\to~\infty
\qquad\hbox{as}\quad \rho_-\to 0\,.\eeq
In other words,  when an infinitesimal 2-wave
crosses a 1-shock of fixed strength $\sigma_1$, its size is amplified by a factor
$\eta'\approx\kappa\,\rho_-^{-2/3}$
which becomes arbitrarily large as
$\rho_-\to 0$, i.e.~as the density approaches vacuum.

As a special case,
if a 1-shock of strength $\sigma_1=1$ crosses a small 2-wave (either compression or rarefaction) of strength $\ve\approx 0$,
the  strength of the shock remains constant, while the strength of the 2-wave is amplified by a factor
\bel{amp}
{\eta(\ve)\over\ve} ~>~{1\over \rho^{2/3}_-}\,.\eeq
The above estimate is valid as soon as the density $\rho_-$ of the left state
(i.e., ahead of the 1-shock) is sufficiently small.
\v
One more case will be of relevance.
Consider a small 1-shock with right state $(u_+, \rho_+) = (0,1)$.
By (\ref{RH2}) its left state $(u_-, \rho_-)$  satisfies
\bel{ss5}
u_- ~=~ \sqrt{ \left(1- {1\over \rho_-}\right)
{\rho_-^3-1\over 3}  }\,.\eeq
Taking $\rho_- = 1-s$, we obtain
\bel{ss6}
u_- ~=~ \sqrt{ \left(s^2 - s^3 + {s^4\over 3}\right)(1+s+s^2 + \cdots)}~=~s\left( 1+ {s^2\over 6}+ {s^3\over 6} + o(s^3)\right)\,.\eeq
Next, consider a small 2-shock, again with right state $(\tilde u_+, \tilde \rho_+) = (0,1)$.
Let $(\tilde u_-, \tilde \rho_-)$ be the left
 state. Taking $\tilde \rho_- = 1+r$, by (\ref{ss5})  we now have
\bel{ss7}
\tilde u_- ~=~ \sqrt{ \left(r^2 + r^3 + {r^4\over 3}\right)(1-r+r^2 - \cdots)}~=~r\left( 1+ {r^2\over 6} -{r^3\over 6} + o(r^3)\right)\,.\eeq
Imposing $u_-= \tilde u_-$ yields
$$r ~=~s + {s^4\over 3} + o(s^4).$$
Referring to Fig.~\ref{f:hyp59},
consider the points
$$A_1= (0,1), \qquad C=(u_-, \rho_-),\qquad D=( u_-, \tilde\rho_-),$$
and let
$$B_2~=~\left(u_- +{\tilde\rho_--\rho_-\over 2}\,,~{\tilde\rho_-+\rho_-\over 2}\right)~=~
\left(2s+  {s^3\over 6}+ {s^4\over 3} + o(s^4)\,,~ 1+  {s^4\over 6} + o(s^4)\right)  $$
Then
the slope of the segment  $A_1 B_2$ is
\bel{slo} {\ds{s^4\over 6} + o(s^4)\over\ds  2s+  {s^3\over 6}+ {s^4\over 3}
+ o(s^4)}~=~
{s^3\over 12}+o(s^3).\eeq

\section{Head-on interactions}
\label{sec:2}
\setcounter{equation}{0}

\begin{figure}[htbp]
\centering
  \includegraphics[scale=0.5]{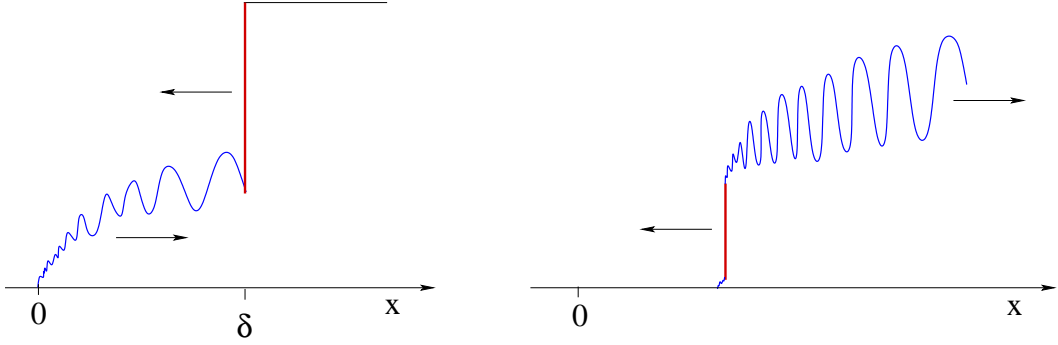}
    \caption{{\small The head-on interaction of a 2-shock with a
    train of smooth 1-waves, in Example~1.}}
\label{f:hyp52}
\end{figure}

{\bf Example 1.}
Consider an initial data $(\bar \rho,\bar u)$ consisting of a 1-shock of strength $\sigma_1=1$,
approaching a train of smooth 2-waves, near vacuum (Fig.~\ref{f:hyp52}).

In terms of Riemann invariants $(w_1, w_2)$, assume that the initial data is given by
\bel{ex1}\left\{ \bega{rl}w_1(0,x)& =~\bar w_1(x)~=~0\,,\cr
w_2(0,x) &=~\bar w_2(x)~=~
x^\alpha\cdot (2+\sin x^{-\beta})\enda \right. \qquad\hbox{for}~~x\in [0,\delta],
\eeq
for suitable constants $\alpha,\beta>0$ and some $x_0>0$ suitably small.
In addition, assume that this  initial data has a
1-shock of size $\sigma_1=1$ at $x=1$, and is constant on the half lines where
 $x<0$ and $x>1$.

For this initial data we construct an approximate solution such that, at every
interaction, the strength of outgoing waves is the same  as in an exact solution.
 However, instead of (\ref{2}) or (\ref{RH5}),
we let the waves travel with  constant speeds, say,
$-c$ for the 1-shock and $c$ for the 2-waves, for some constant $c>0$.
Then, as $t\mapsto T=\delta/2c$, all the 2-waves  cross the shock
and the total variation of the
solution approaches infinity.   Indeed, by (\ref{amp}), choosing
$\delta>0$ sufficiently small  the following holds.  At any  time $0<t<T$
we have
\bel{ampx}
\bega{rl}|w_{2,x} (t,x+ct) |&\ds\geq~{1\over
\bar \rho(x)^{2/3}} \cdot |\bar w_{2,x}(x)|
\cr\cr
&\geq~\ds
{1\over \left(3x^\alpha/2\right)^{2/3}} \cdot \Big| \alpha x^{\alpha-1} (2+\sin x^{-\beta})
-\beta x^{\alpha-\beta-1}\cos x^{-\beta}\Big|,\enda\eeq
for all $ x\in[0, x_0]$ such that $x+ct>1-ct$.

If we now choose $0<\alpha/3<\beta<\alpha<1$, then the initial data has bounded variation, because
$$\int_0^1 |w_{2,x} (x) | \,dx~=~\int_0^1\Big| \alpha x^{\alpha-1} (2+\sin x^{-\beta})
-\beta x^{\alpha-\beta-1}\cos x^{-\beta}\Big|\, dx~<~\infty.$$
On the other hand, as $t\to T-$ the total variation blows up, because
$$\bega{l}\ds\lim_{t\to 1/2c}\int_{1-2ct}^{x_0} |w_{2,x} (t, x+ct) | \,dx\cr\cr
\qquad \ds=~
\left({2\over 3}\right)^{2\over 3}\int_0^{x_0} \Big| \alpha x^{(\alpha/3)-1} (2+\sin x^{-\beta})
-\beta x^{(\alpha/3)-\beta-1}\cos x^{-\beta}\Big|\, dx~=~+\infty.\enda
$$
\v
In the previous example the initial data contains vacuum.
Moreover, in the terminal
profile the large variation is achieved at very low gas density.
The next example shows that this arbitrarily large amplification
of the total variation can be achieved also with an initial datum
having uniformly positive density. As before, we require that at each interaction
the strength of outgoing waves is the same as in an exact solution, but we allow
a small error in the wave speed.   Namely, in our approximate solution
all 1-waves travel with speed $-c$ while 2-waves travel with speed $c$.
\v
\begin{figure}[htbp]
\centering
  \includegraphics[scale=0.4]{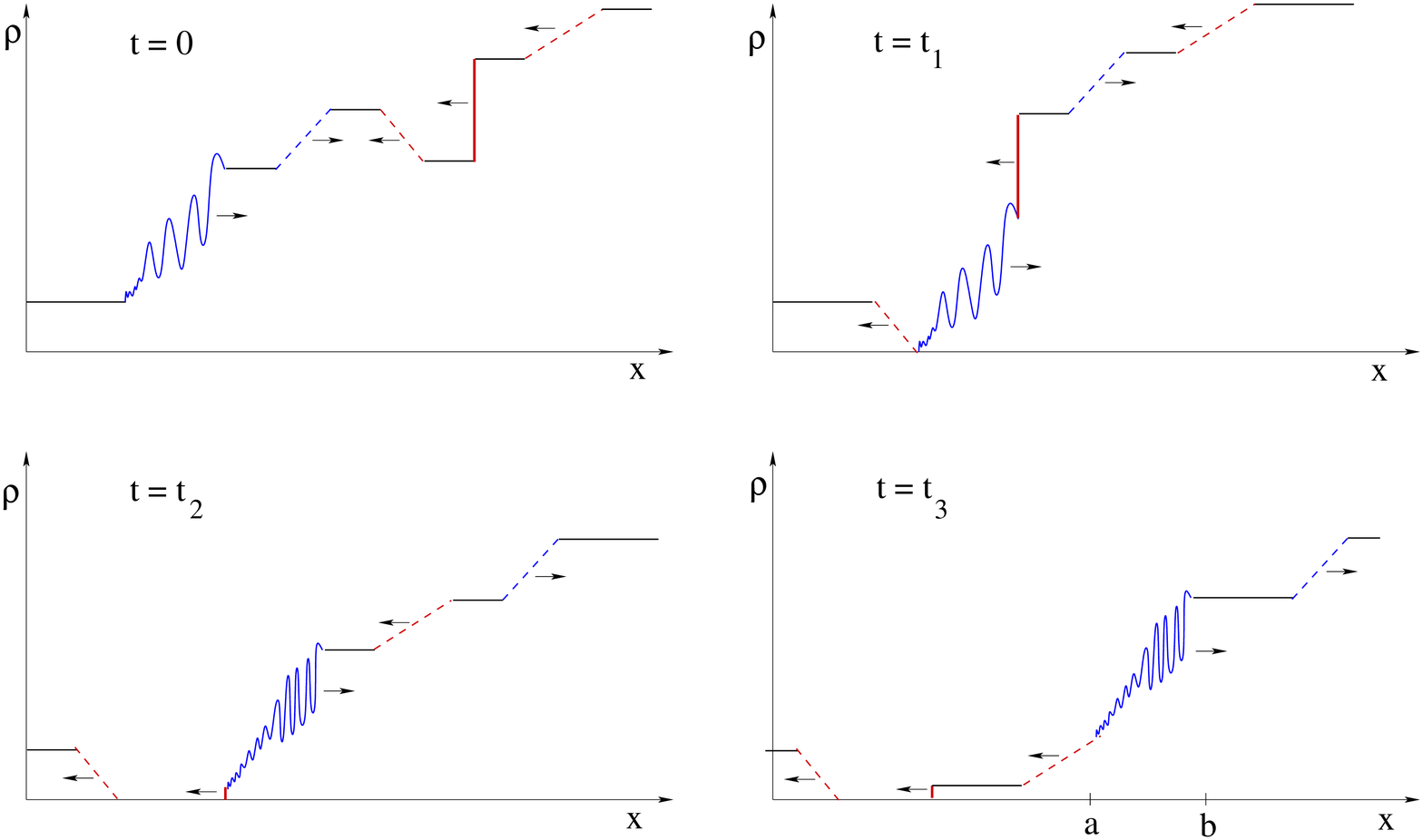}
    \caption{{\small Evolution of the density profile in Example 2.
    Top left: at time $t=0$ the initial density is uniformly positive. Top right:
    after two rarefaction waves of the opposite families cross, in the middle section
    the same configuration in Example 1 is recreated. Bottom left: when the 1-shock
    crosses the train of 2-waves at low density, the total variation grows without bounds. Bottom right: after crossing a 1-compression, an arbitrarily large
    amount of total variation occurs within the interval $[a,b]$, at uniformly positive density.  }}
\label{f:hyp87}
\end{figure}

\begin{figure}[htbp]
\centering
  \includegraphics[scale=0.5]{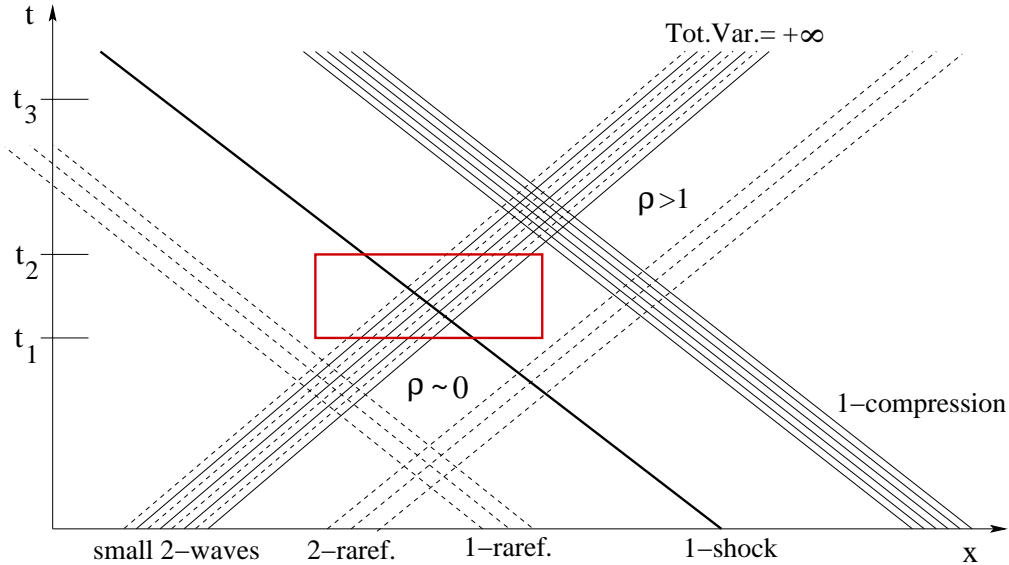}
    \caption{{\small The wave pattern described in Example 2.
    The box shows the region where the interactions in
Example 1 take place.}}
\label{f:hyp46}
\end{figure}

{\bf Example 2.}  As shown in Fig.~\ref{f:hyp87},
consider an initial data similar to (\ref{ex1}),
but with the insertion of two additional rarefaction waves, and a compression.
At time $t_1$, after crossing rarefaction waves of the opposite family, at time
$t=t_1$ the train of 2-waves and the 1-shock recreate the initial data in Example~1.
At time $t_2$, the total variation of the second Riemann invariant $w_2(t,\cdot)$
becomes infinite.    At time $t_3$, after crossing a 1-compression,
this infinite total variation occurs at uniformly positive density.
In the $t$-$x$ plane, the solution is described in Fig.~\ref{f:hyp46}.

\v
{\bf Remark 1.} Example 2 shows that there exists an initial profile
$(\bar u, \bar \rho)\in BV$ with $\bar \rho(x)\ge\rho_0>0$,
and an approximate solution with fronts moving at constant speed $\pm c$
such that the following holds.  If at each interaction the strength of
outgoing waves
is the same as in the exact solution of the Riemann problem,
then the total variation blows up in finite time.  More precisely, at some time $t_3$
\bel{ITV}
\rho(t_3,x)\geq 1\qquad \hbox{for}~~x\geq 0,\qquad\quad \TV\Big\{ \rho(t_3,\cdot)\,;
~~[0,\infty[\,\Big\}~=~\infty.\eeq
In particular, is is not possible to put a continuous
weight on wave strengths, possibly approaching infinity as $\rho\to 0$,
in order to control the total variation.  Indeed, for any choice of the weights,
continuous on the region where $\rho>0$, the initial weighted total variation
will be finite, while the final weighted total variation will be infinite.

\v

\v

\section{An example with uniformly positive density}
\label{sec:3}
\setcounter{equation}{0}

\begin{figure}[htbp]
\centering
  \includegraphics[scale=0.4]{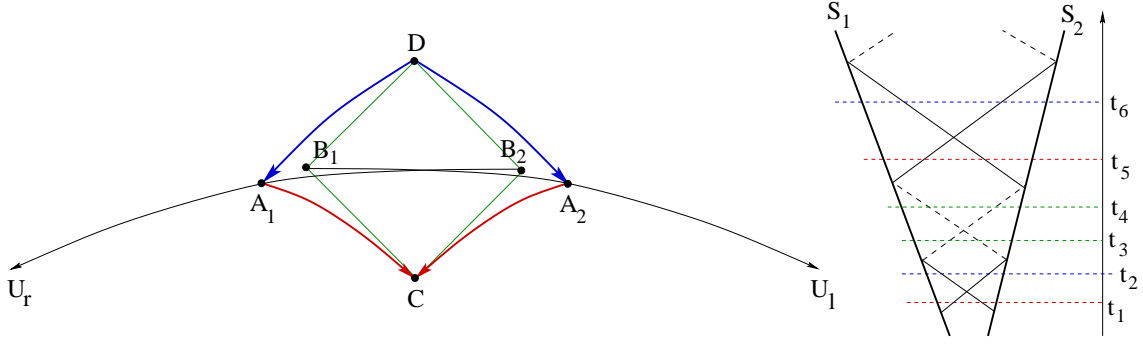}
    \caption{{\small A periodic interaction pattern.}}
\label{f:hyp59}
\end{figure}

In the previous examples, the blow up of the total variation
was achieved because
waves crossing a shock of unit strength were amplified by an arbitrarily large factor
as the gas density approached vacuum.
The following example shows that the total variation can become
arbitrarily large even  if the gas density remains uniformly bounded away from zero,
at all times.
\v
{\bf Example 3.}   STEP 1.
We begin by constructing
a symmetric interaction pattern containing four wave fronts,
as shown in Fig.~\ref{f:hyp59}.   We choose the strengths of the two
large shocks $S_1, S_2$ and of the two intermediate waves in such a way that,
after a whole round of interactions, these strengths are the same as at the initial time.  Working in the $u$-$\rho$ plane, this is done as follows.

\begi
\item[(i)] Start by constructing two symmetric shocks: the 1-shock
$A_1C$ and the 2-shock $A_2C$, approaching each other.
\item[(ii)] Determine the two outgoing shocks $DA_1$ and $DA_2$, resulting
from the crossing of the above two shocks.
\item[(iii)] Construct a square  having two opposite vertices
at $C$ and $D$.   Call $B_1$, $B_2$ the remaining two vertices.
\item[(iv)] Choose  $U_l$ so that the two points $B_1$ and $A_2$
are on the same 1-shock curve with left state $U_l$.
Symmetrically, choose  $U_r$ so that the two points $B_2$ and $A_1$
are on the same 2-shock curve with  right state $U_r$.
\endi

\begin{figure}[htbp]
\centering
  \includegraphics[scale=0.5]{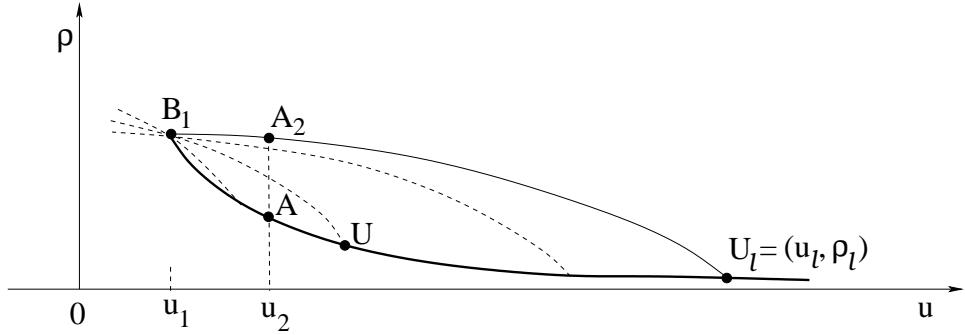}
    \caption{{\small By moving the point $U$ along the 1-shock
    curve with right state $B_1$, we eventually
   reach a left state $U_l$     such that the 1-shock curve through $U_l$
   contains $A_2$ as well.}}
\label{f:hyp62}
\end{figure}

The existence of  states $U_l$, $U_r$ satisfying the conditions (iv) is now proved
(see  Fig.~\ref{f:hyp62}).

{\bf Lemma 1.} {\it In the $u$-$\rho$ plane, consider two points $B_1=(u_1,\rho_1)$
and $A_2=(u_2, \rho_2)$. Assume that
\begi
\item[(i)]~~$u_1<u_2$, and $\rho_1>\rho_2$.
\item[(ii)] Calling $A=(u_2, \rho_2^*)$ the point on the 1-shock curve with right state
 $B_1$
with the same $u$-component as $A_2$, one has $\rho_2^*<\rho_2$.
\endi
Then there exists
a unique $U_l= (u_l,\rho_l)$, with $0<\rho_l<\rho_2$, such that both $B_1$ and $A_2$ lie on the 1-shock curve
with left state state $U_l$.}

{\bf Proof.}
We shall use (\ref{RH2}) with
$(u_-,\rho_-) = (u_l, \rho_l)$ while $(u_+, \rho_+) = (u_1, \rho_1) $ or $(u_+, \rho_+) = (u_2, \rho_2) $.
To prove the lemma we need to find $(u_l, \rho_l)$
such that
\bel{url}
u_l-u_1~=~\sqrt{ \left({1\over\rho_1}- {1\over \rho_l}\right)\left(
{\rho_l^3\over 3}-{\rho_1^3\over 3}\right)}\,,
\qquad
u_l-u_2~=~\sqrt{ \left({1\over\rho_2}- {1\over \rho_l}\right)\left(
{\rho_l^3\over 3}-{\rho_2^3\over 3}\right)}\,.\eeq
Equivalently,
\bel{url2}
u_2-u_1~=~G(\rho_l)~\doteq~\sqrt{ \left({1\over\rho_1}- {1\over \rho_l}\right)\left(
{\rho_l^3\over 3}-{\rho_1^3\over 3}\right)}-
\sqrt{ \left({1\over\rho_2}- {1\over \rho_l}\right)\left(
{\rho_l^3\over 3}-{\rho_2^3\over 3}\right)}.\eeq
The assumption (ii) implies
$$G(\rho_2) ~=~ \sqrt{ \left({1\over\rho_1}- {1\over \rho_l}\right)\left(
{\rho_l^3\over 3}-{\rho_1^3\over 3}\right)}~<~u_2-u_1\,.$$
Moreover, a direct computation yields
$${\partial\over\partial \rho_l} G(\rho_l) ~<~0,
\qquad \lim_{\rho_l\to 0} G(\rho_l) ~=~+\infty.$$
Therefore there exists a unique value of $\rho_l$ such that $G(\rho_l) = u_2- u_1$.
\endproof

Having constructed the above wave curves,
consider the following interaction pattern, shown in Fig.~\ref{f:hyp59}, right:
\begi
\item At time $t_1$ the initial datum consisting of four shocks, connecting the states
$U_l, A_2, C, A_1, U_r$.
\item At time $t_2$ the profile still consists of four shocks
(the two middle ones
have crossed each other), connecting the states
$U_l, A_2, D, A_1, U_r$.
\item At time $t_3$ the profile consists of two large shocks and two rarefactions
(the two middle shocks have joined the big ones, generating two rarefactions),
connecting the states
$U_l, B_1, D, B_2, U_r$.
\item At time $t_4$ the profile still consists of two large shocks and two rarefactions
(the two rarefactions have crossed each other),
connecting the states
$U_l, B_1, C, B_2, U_r$.
\item At time $t_5$ the initial datum consisting of four shocks
(the two rarefactions have impinged on the big shocks, generating two
intermediate shocks of the opposite families), connecting the states
$U_l, A_2, C, A_1, U_r$, exactly the same as at time $t_1$.
The pattern thus repeats itself.
\endi
\v
\v
{\bf Remark 3.} It is clear that for the above example it is essential to have large
total variation.   Indeed, if we choose the middle shock $A_1C$
small, then the line $A_1 B_2$ will be almost horizontal and the point $U_r$
must be far away, with density close to zero.   \\
On the other hand, in this solution obtained by front tracking
the density trivially remains uniformly bounded away from zero.
Since the interaction pattern is periodic, we conclude that
\v
{\it
Even under the assumption that the density remains uniformly positive,
there is no way to construct a Lyapunov functional controlling the total variation,
which is strictly decreasing at every interaction.}

In the above approximate solution the wave strengths follow a periodic pattern.
To achieve an arbitrarily large amplification of the total variation, a further
construction is needed.

\begin{figure}[htbp]
\centering
 \includegraphics[scale=0.5]{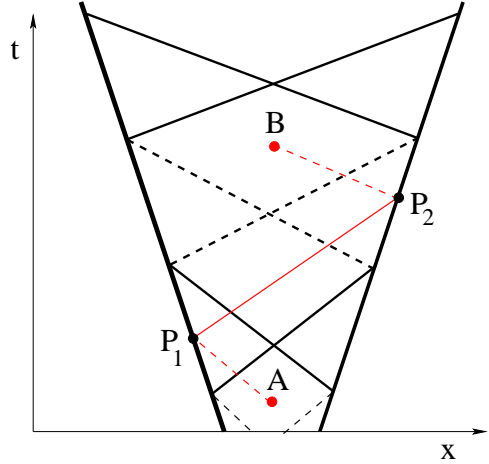}
    \caption{{\small A periodic pattern that amplifies a small wave front.}}
\label{f:hyp88}
\end{figure}

\v
STEP 2:  As shown in Fig.~\ref{f:hyp88}, right, on top of the previous pattern
we add a very small wave front.  To fix the ideas, consider a 1-rarefaction
of strength $\ve>0$, located at $A$.     Within a time period, this front will
\begi
\item Cross the intermediate 2-shock.
\item Interact with the large 1-shock at $P_1$ producing a 2-compression.
\item Cross the intermediate 1-shock.
\item Cross the intermediate 1-rarefaction.
\item Interact with the large 2-shock at $P_2$ producing a 2-rarefaction.
\item Cross the intermediate 2-rarefaction.
\endi
We analyze the case where the two middle shocks are small and
the density of their left state is $\approx 1$.
We claim that, when the additional front reaches $B$, its size will be increased
 by a factor $\kappa>1$.

Indeed, when a small wave of strength $\ve^-$ crosses a shock of the opposite family
of strength $s$ at density $\rho_-=1$, by (\ref{epex}) the strength of the outgoing front is
\bel{amp1}\ve^+ ~=~\left(1+ {s^3\over 3}+o(s^3)\right)\ve^-.\eeq

When the front crosses a rarefaction of the opposite family, its strength  does
not change.

Finally,  when the  small wave impinges on a large shock at $P_1$ or at $P_2$,
we need to estimate  the relative size of the reflected wave front.
Toward this goal, let
$\rho = \psi(u)$ be the equation of the shock curve with right state $U_r$,
passing through both $A_1$ and $B_2$, as constructed in Lemma 1.
Calling $\ve^-$, $\ve^+$ the strengths of the front before and after interaction,
to leading order we have
\bel{reflec}\ve^+~=~(1-2\psi'(u)) \ve^-~=~\left(1- {s^3\over 6}+o(s^3)\right)\ve^-.\eeq
Indeed, by  (\ref{slo}) we have $\psi'(u)= s^3/12 + o(s^3)$.
Calling $\ve_A$ and $\ve_B$ respectively the strengths of  the small wave-front
at $A$ and at $B$, we thus have
\bel{ampe}
\ve_B~=~\left(1+ {s^3\over 3}+o(s^3)\right)^2\left(1- {s^3\over 6}+o(s^3)\right)^2
\ve_A~=~\left(1+ {s^3\over 3}+o(s^3)\right)
\ve_A\,.\eeq

\begin{figure}[htbp]
\centering
  \includegraphics[scale=0.5]{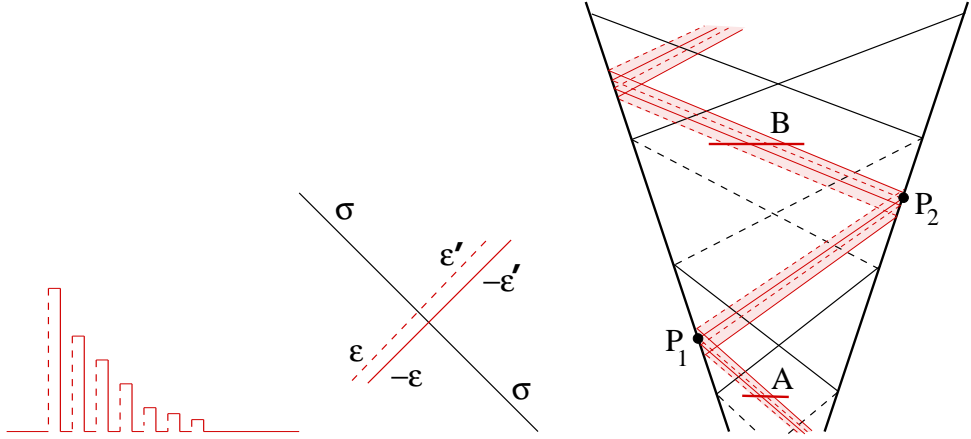}
    \caption{{\small
 Right: a periodic pattern that amplifies infinitesimal waves.
 Center: if a front of strength $\sigma$ crosses  a compression (not a shock!) of size $-\ve$ and then a rarefaction of size $\ve$,
its final strength is not changed.
 The outgoing fronts have strengths $\ve', -\ve'$.  Here $\ve'>\ve$ if the
 front $\sigma$ is a shock, otherwise $\ve'=\ve$.
 Left: to construct a interaction pattern that yields an arbitrarily large
 total variation, we replace the single infinitesimally small front
 in Fig.~\ref{f:hyp88}  by
 countably many pairs rarefaction + compression, of opposite size. The total strength
 of these waves is finite, each front having strength $\leq\ve$.}}
\label{f:hyp85}
\end{figure}

\v
STEP 3.  Consider a periodic pattern that amplifies an infinitesimal
wave, as in Step 2.  By continuity, there exists $\lambda>1$ and $\ve>0$
such that any wave-front (rarefaction or compression)
of size $<\ve$, traveling from $A$ to $B$ along the path
in Fig.~\ref{f:hyp88},  is amplified by a factor $\geq\lambda$.

We now construct an initial set of wave fronts where the infinitesimal  front
is replaced by countably many pairs ``rarefaction + compression", whose sizes
exactly cancel each other (Fig.~\ref{f:hyp85}, left).
At time $t=0$, the total strength of all these small
fronts can be taken to be $=1$.
The key observation is that each of these pairs {\it leaves no footprint} on the
underlying solution constructed in Step 1.   Indeed, if a front of size
$\sigma$  crosses
a rarefaction and a compression of exactly opposite sizes, after the two crossings the size of the front is still $\sigma$ (Fig.~\ref{f:hyp85}, center).
As a result, the pattern of four large fronts retains its periodicity.

By construction, after each period each pair of opposite small wavefronts
is enlarged by a factor $\geq\lambda$.   When a pair grows to size $>\ve$, we can
perform a partial cancellation so that its size remains $\in [\ve/2,~\ve]$.

Since the total number of small wave-fronts is infinite,
after several periods a larger and larger number of  pairs (compression + rarefaction)
reaches size $>\ve/2$.   Hence, as $t\to\infty$,
 the total variation of this approximate solution grows without bounds.

\section{Blow up in finite time}
\label{sec:4}
\setcounter{equation}{0}

The previous construction shows that, if the total variation is initially sufficiently
large,
then there exists an interaction pattern that renders
the total variation arbitrarily large as $t\to\infty$.  This can be achieved even
with a uniform lower bound on the density.

The next question is whether one can arrange the order of wave-front interactions
so that the total variation blows up in finite time.   Notice that this is not the case
in the previous example.  Indeed, if 1-fronts travel with speed $\dot x\in [-C, 0[\,$
and 2-fronts travel with speed $\dot x_2\in \,]0,C]$, it takes a uniformly positive
amount of time for each intermediate front to
bounce back and forth between two large shocks.  Hence the arbitrarily large amplification
of the total variation is only achieved in the limit as $t\to +\infty$.

\begin{figure}[htbp]
\centering
  \includegraphics[scale=0.5]{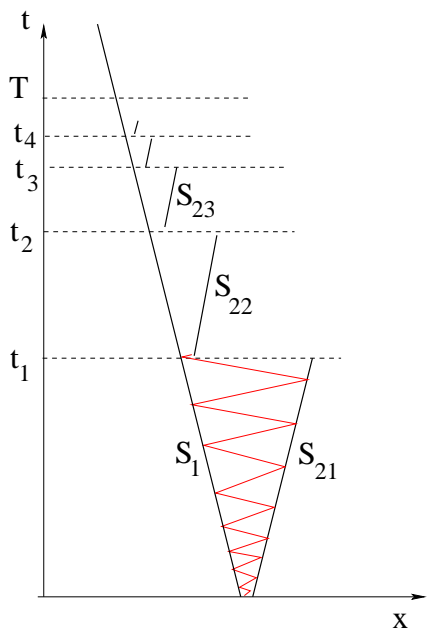}
    \caption{{\small A pattern yielding finite time blow up of the total variation.
    At each step the total amount of small waves bouncing back and forth between the two large shocks keeps increasing.
}}
\label{f:hyp84}
\end{figure}

\begin{figure}[htbp]
\centering
  \includegraphics[scale=0.4]{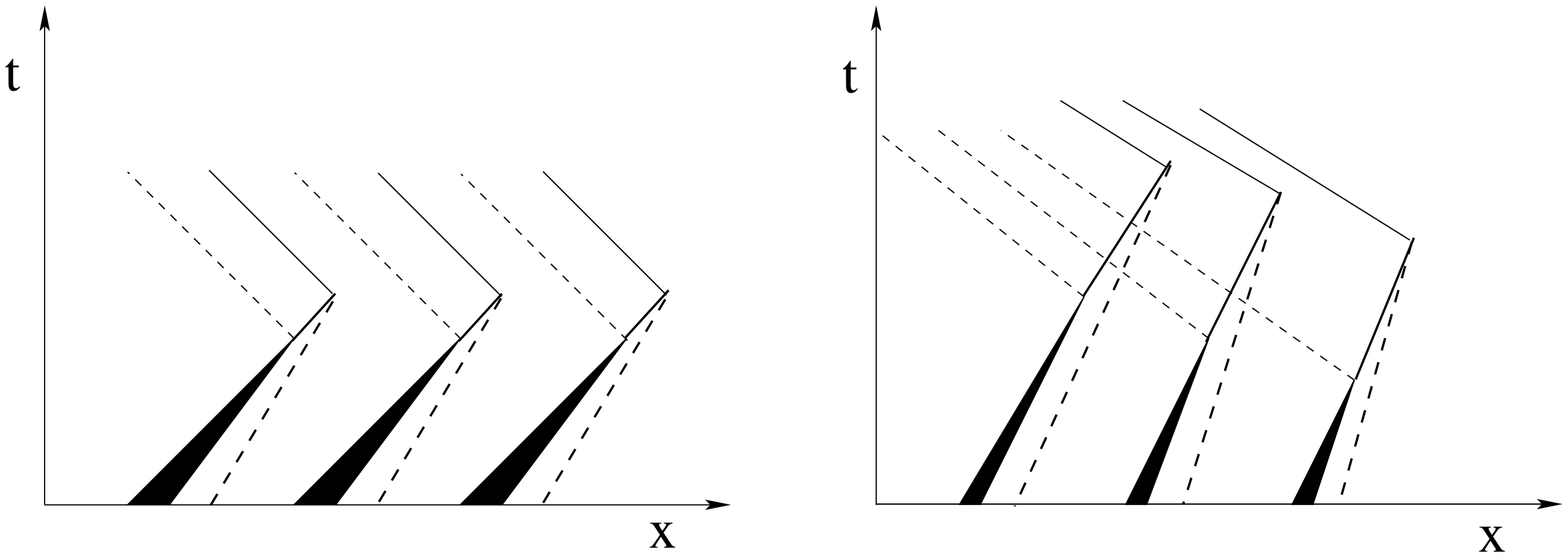}
    \caption{{\small Starting with a large number of pairs of small waves,
    one can generate a large number of similar pairs in the other characteristic family
    (left), or one single large pair of fronts (right).
      }}
\label{f:hyp89}
\end{figure}

\begin{figure}[htbp]
\centering
  \includegraphics[scale=0.35]{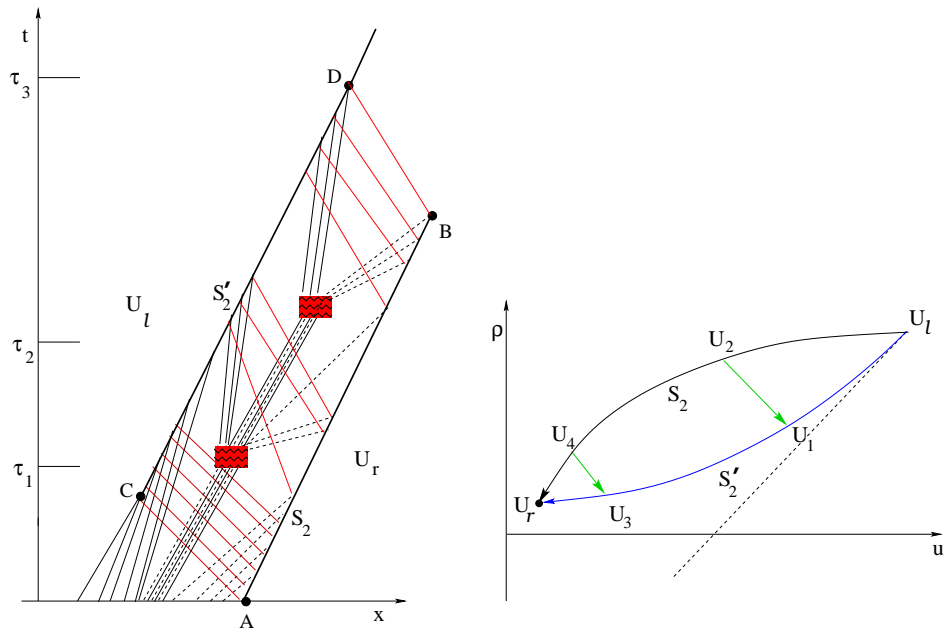}
    \caption{{\small Left:
    between $A$ and $B$ the shock $S_2$ is completely canceled by impinging
    2-rarefactions.   Between $C$ and $D$, the
    new 2-shock $S_2'$ is formed by impinging compression waves.   Notice that
    the 1-compression fronts emerging from $S_2$ are used to completely cancel
    the 1-rarefactions that would otherwise be produced by the interactions
    of 2-compressions with $S_2'$.   Right:  if the shock $S_2$ is very large,
    it cannot be canceled by one single large rarefaction front.   Therefore,
    we need to produce several small rarefactions at subsequent times,
    so that the density $\rho$ remains  uniformly positive.
      }}
\label{f:hyp86}
\end{figure}

In this section we briefly indicate how the previous construction can be modified,
providing finite time blow up of the total variation.  The main idea is illustrated in
Fig.~\ref{f:hyp84}.    We consider a sequence of times $0=t_0<t_1<t_2<\cdots<T$.
During each time interval $J_i = [t_{i-1}, t_{i}]$, a countable number of pairs
of small waves (compression + rarefaction) is amplified by a very large factor.
Before time $t_{i}$, the large 2-shock  $S_{2i}$ is completely canceled
by impinging 2-rarefactions, and a new 2-shock $S_{2(i+1)}$ of the same strength
is recreated at a location closer to the large 1-shock $S_1$.
Figures~\ref{f:hyp89} and  \ref{f:hyp86}
show how this can be achieved, starting with very many
pairs of small waves (compression + rarefaction).
By letting each compression front collapse to a shock, and then canceling
this shock with a rarefaction front of the same family,
we obtain a train of pairs of small waves (compression + rarefaction)
in the opposite family (Fig.~\ref{f:hyp89}, left).
By varying  the locations of these interactions,
instead of many pairs of small waves we can achieve
a large compression
followed by a large rarefaction front (Fig.~\ref{f:hyp89},
right).

The basic step is illustrated in
Fig.~\ref{f:hyp86}.  A large number of small compression+rarefaction
pairs produces
    a large 2-rarefaction, which starts depleting the 2-shock along the line
    $AB$, and a large compression, which builds up a new 2-shock along the line
    $CD$.
    Since we cannot allow the density to become negative, it may not be possible to
    cancel the large shock $S_2$ with one single large 2-rarefaction.
    For this reason, this cancellation may be accomplished in several stages.

  For example,   the first set of 2-rarefactions reduce the
    size of the 2-shock $S_2$ from  $(U_\ell, U_r)$ to $(U_2,U_r)$.
    At a time $\tau_1$, the profile $u(\tau_1,\cdot)$ thus contains the 2-shock
    $(U_\ell, U_2)$, the 1-compression $(U_2, U_1)$, and the new 2-shock
    $(U_1, U_r)$.

    At a later time $\tau_2$, the profile $u(\tau_2,\cdot)$  contains the (shrinking)
    2-shock
    $(U_\ell, U_4)$, the 1-compression $(U_4, U_3)$, and the (growing) 2-shock
    $(U_3, U_r)$.

    At a later time $\tau_3$, the original 2-shock has been completely depleted by
    2-rarefactions.   A 2-shock connecting exactly the same two states
    $(U_\ell, U_r)$ is formed at a different location, as desired.

By canceling the  large 2-shock and reconstructing it at a different location, shifted to the left, we can reproduce the pattern in Fig.~\ref{f:hyp84}.
Since at each step the total strength of the small intermediate waves can be amplified
by an arbitrarily large factor, as $t\to T$ the total variation of our
approximate solution blows up to $+\infty$.


\section{Concluding remarks}
\setcounter{equation}{0}
The examples presented in this paper show that, if the strength of wave-fronts
is computed exactly but some error is allowed their speeds, then
the total variation of approximate solutions can blow up. It is interesting to
revisit some of the previous examples, taking into account the decay of
rarefaction waves due to genuine nonlinearity.
Looking at exact solutions, it becomes clear that these particular 
interaction patterns do not yield an arbitrarily large
amplification of the total variation.

\subsection{Head-on interactions, near vacuum.}
Consider an exact solution of the system (\ref{1}), with initial
data as in Example 1. We show that there is no way to choose
$\alpha,\beta$ in (\ref{ex1}) so that the following requirements are simultaneously
satisfied:
\begi
\item[(i)] The 1-shock  crosses all 2-waves in finite time.
\item[(ii)] The 2-waves do not break before crossing the shock.
\item[(iii)] The sum of strengths of the 2-waves is initially finite, and becomes infinite as they all cross the 1-shock.
\item[(iv)] The 2-waves do not break immediately after crossing the shock.
\endi

For a 1-shock of unit strength, assume that the left state (ahead of the shock)  has
density $\rho_-\approx 0$.
 Then by (\ref{RH2}) the right state (behind the shock) has density
$\rho_+$ satisfying
$${\rho_+^3\over 3}\cdot {1\over \rho_-} ~\approx 1.$$
Hence the right state and the speed of the shock are given 
respectively by
\bel{r-}
\rho_+~\approx~(3\rho_-)^{1/3}\,,\qquad \dot x~\approx~-{1\over \sqrt 3} \sqrt{\rho_+^3
 \rho_-}~=~3^{-1/3}\, \rho_-^{2/3}.\eeq
If the initial profile is given by (\ref{ex1}), so that
\bel{iprof}
w_1(0,x)~=~ 0,\qquad w_2(0,x)~=~x^\alpha(2+\sin x^{-\beta})
\qquad \qquad x\in [0,\delta],\eeq
then  the requirement (i) will be satisfied
provided that 
\bel{r1} {2\over 3}\alpha~<~1\,.\eeq
Next, to make sure that the 2-waves do not break before crossing the 1-shock,
we look at the evolution of $w_{2,x}$ along characteristics. From
$$w_{2,t} + \left({w_2\over 2}\right)^2 w_{2,x}~=~0$$
it follows
$$w_{2,xt} + \left({w_2\over 2}\right)^2 w_{2,xx}~=~- {w_2\over 4}\, w_{2,x}^2\,.$$
By a comparison argument,
we conclude that the gradient $w_{2,x}$ will not blow up before time $T>0$ 
provided that the initial data satisfy
$$|w_{2,x}(0,x)| \cdot w_2(0,x)~<~T^{-1}.$$
Recalling (\ref{iprof}),  the condition (ii) is thus satisfied if
\bel{r2}
\alpha + (\alpha-\beta-1)~>~0\,.\eeq
As shown in the discussion of Example 1, condition (iii) is satisfied
provided that 
\bel{r3} 0~<~{\alpha\over 3}~<~\beta~<~\alpha~<~1.\eeq

To check whether (iv)
can be satisfied,
let $T$ be the time when the 1-shock reaches the origin, crossing all
2-waves.    
Denote by $t\mapsto x(t,y)$ the position of a 2-characteristic starting at 
$x(0,y)=y$,
Calling $\tau(y)$ the time where this 2-characteristic 
crosses the 1-shock, we find
$$T-\tau(y)~\approx~y^{1-2\alpha/3}.$$
We consider the evolution of $w_{2,x}(t,x(t,y))$ along this
 2-characteristic. 
 For $t<\tau(y)$ we have 
$$|w_{2,x}(t,x(t,y))|~\approx~y^{\alpha-\beta-1}. $$
When this characteristic crosses the 1-shock at time $\tau=\tau(y)$,
by (\ref{amp}), this gradient is
amplified by a factor 
$${w_{2,x}(\tau+)\over
w_{2,x}(\tau-)}~ \approx~ y^{-2\alpha/3}.$$
Moreover, $w_2(\tau+)\approx y^{\alpha/3}$.
To make sure that this gradient remains bounded during the time interval $[\tau(y), T]$,
  we need
$$y^{-2\alpha/3}\cdot |w_{2,x}(\tau-)|\cdot  w_2(\tau+)\,(T-t)~=~\O(1).$$
Therefore we should have
\bel{r4}
y^{-2\alpha/3} \cdot y^{\alpha-\beta-1} \cdot y^{\alpha/3}\cdot y^{1-2\alpha/3}~
=~y^{- \beta} ~=~\O(1).\eeq
This condition is incompatible with the requirement (\ref{r3}) that $\beta>0$.


\begin{figure}[htbp]
\centering
  \includegraphics[scale=0.39]{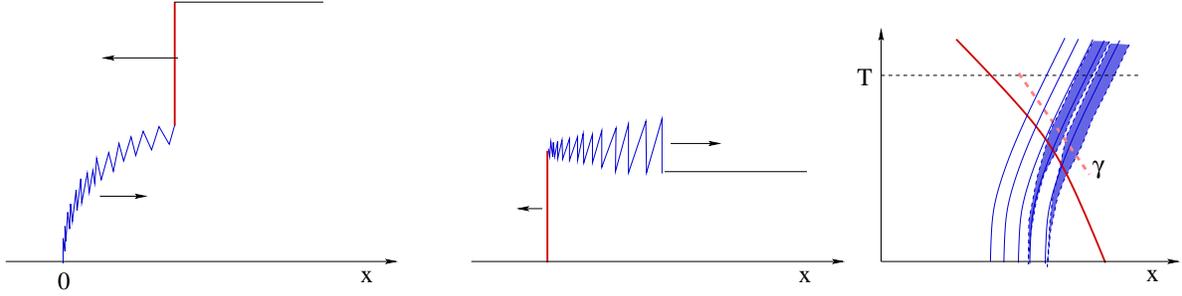}
    \caption{{\small Left: Before crossing the large 1-shock, the small 2-waves 
    do not break, because when the gas density $\rho>0$
    is very small, the system is almost linearly degenerate.
    Center: after crossing the 1-shock, the 2-compression waves break and 
 a large amount of cancellation between 2-rarefactions and 2-shocks occurs.    At the terminal time $T$ the total strength of waves is still finite, due to these cancellations.  Right:  in the $t$-$x$ plane this pattern  produces an
    infinite total variation only along the (dashed) time-like curve $\gamma$.  }}
\label{f:hyp77}
\end{figure}

\subsection{Waves bouncing back and forth between two large shocks}

\begin{figure}[htbp]
\centering
  \includegraphics[scale=0.45]{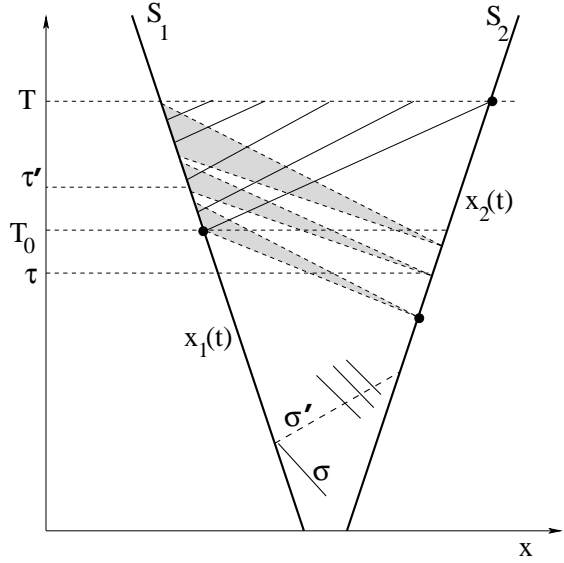}
    \caption{{\small  The total amount of compression waves (including shocks)
    at time $T$ is bounded in terms of the total amount of rarefaction
    waves that impinge on the large 1-shock during the interval $[T_0,T]$.}}
\label{f:hyp78}
\end{figure}

From our earlier analysis, this should be 
the pattern that achieves the greatest amplification
of wave strengths (Fig.~\ref{f:hyp78}).  If  the size of the shocks $S_1, S_2$
is sufficiently large, the strength of a reflected 2-front $\sigma'$ is almost the same 
as the strength of the impinging 1-front $\sigma$.    Afterwards, as this  2-front
crosses other 1-shocks, its strength increases by a large factor.
Repeating this pattern, it may appear that an arbitrarily large amplification of wave strengths can be achieved.  The following analysis shows that this is not the case,
if we take into account the decay of rarefaction waves due 
to genuine nonlinearity.

For some constant $c_0$, the two large shocks will have speeds
\bel{sspd}
\dot x_1(t)~\le~-c_0 ~<~0~<~c_0~\leq~ \dot x_2(t)\,.\eeq
Consider a 1-rarefaction wave (Fig.~\ref{f:hyp78}) emerging from
the large 2-shock at some time $\tau$ and impinging on the
opposite 1-shock at time $\tau'$.    The upper and lower estimates on the
velocity yield an estimate of the form
\bel{time}
\tau'~\geq~\kappa \tau\,.\eeq
By wave decay estimates, the density of such 1-rarefaction at time $\tau'$
is $$\leq~ C\cdot {1\over \tau'-\tau}~\le~{C\kappa\over \kappa-1}\,
{1\over\tau'}\,.$$
Therefore, the total amount of 1-rarefactions that impinge on the large 1-shock
within a time interval $[T_0,T]$ is
$$\le~\int_{T_0}^{T} {C\kappa\over \kappa-1}\,{1\over\tau'}\, d\tau'.$$
An entirely similar estimate holds for the 2-rarefactions impinging on the large
2-shock.

Next, fix any large time $T$.  As in Figure~\ref{f:hyp78},
consider the maximal backward 1-characteristic
through the point $(T, x_2(T))$.   This will cross the large 1-shock at an earlier time
$T_0$.
The total amount of
2-shocks (together with 2-compressions) at a given time $T$ can be estimated as:
\bel{ests}\bega{c}
\hbox{Total amount of
2-shocks
at time $T$ contained in the interval $[x_1(T), x_2(T)]$}\cr\cr
=~\O(1)\cdot \hbox{amount of 1-rarefactions
impinging on the  1-shock for $t\in [T_0,\,T]$}\cr\cr
\ds=~\O(1)\cdot~\int_{T_0}^{T} {C\kappa\over \kappa-1}\,{1\over\tau'}\, d\tau'
~=~\O(1)\cdot \ln{T\over T_0}~=~\O(1).
\enda
\eeq
Here $\O(1)$  denotes a quantity that remans uniformly bounded 
(provided that some upper and lower bounds on the density $\rho$ are given).

An entirely similar estimate of course holds for 1-shocks
(together with 1-compressions).
In addition, the rarefaction waves can be estimated as
\bel{estr}\bega{c}
\hbox{total amount of
2-rarefactions at time $T$ contained in the interval $[x_1(T), x_2(T)]$}\cr\cr
=~\O(1)\cdot \hbox{amount of 1-shocks
impinging on the  1-shock for $t\in [T_0,\,T]$}\cr\cr
\leq~
\hbox{total amount of
1-shocks at time $T_0$ contained in the interval $[x_1(T_0), x_2(T_0)]$}~\leq~\hbox{const.}\cr\cr
\enda
\eeq
Here the last inequality follows from (\ref{ests}), with $T$ replaced by $T_0$.
This yields a uniform a priori bound on the total strength  of waves produced by 
this particular wave-interaction pattern.
\v

\end{document}